\documentclass[number,citesort,dvips]{arxbj}
\usepackage{upgreek}
\usepackage{graphicx}

\aid{0}
\volume{17}
\issue{1}
\pubyear{2011}
\firstpage{170}
\lastpage{193}
\doi{10.3150/10-BEJ265}

\makeatletter

\newtheorem{theo}{Theorem}
\newtheorem{lem}{Lemma}

\newremark{Rem}{Remarks}

\makeatother

\begin{document}
\begin{frontmatter}

\title{Rice formulae and Gaussian waves}
\runtitle{Rice formulae and Gaussian waves}

\begin{aug}
\author[a]{\fnms{Jean-Marc} \snm{Aza\"{\i}s}\thanksref{a}\corref{}\ead[label=e1]{azais@cict.fr}},
\author[b]{\fnms{Jos\'e R.} \snm{Le\'on}\thanksref{b}\ead[label=e2]{jose.leon@ciens.ucv.ve}}
\and
\author[c]{\fnms{Mario} \snm{Wschebor}\thanksref{c}\ead[label=e3]{wschebor@cmat.edu.u}}
\runauthor{J.-M. Aza\"{\i}s, J.R. Le\'on and M. Wschebor}
\address[a]{Universit\'e de Toulouse, IMT, LSP,
F31062 Toulouse Cedex 9, France. \printead{e1}}
\address[b]{Escuela de Matem\'atica, Facultad de Ciencias,
Universidad Central de Venezuela, A.P. 47197, Los Chaguaramos,
Caracas 1041-A, Venezuela. \printead{e2}}
\address[c]{Centro de Matem\'{a}tica, Facultad de
Ciencias, Universidad de la Rep\'{u}blica, Calle Igu\'{a} 4225,
11400 Montevideo, Uruguay. \printead{e3}}
\end{aug}

\received{\smonth{4} \syear{2009}}
\revised{\smonth{1} \syear{2010}}

%
\begin{abstract}
We use Rice formulae in order to compute the
moments of some level functionals which are linked to problems in
oceanography and optics: the number of
specular points in one and two dimensions,
the distribution of the normal angle of level curves and the number of
dislocations in random wavefronts. We compute
expectations and, in some cases, also second moments of such
functionals. Moments of order greater than one are more involved,
but one needs them whenever one wants to perform statistical
inference on some parameters in the model or to test the model
itself. In some cases, we are able to use these computations to obtain
a central limit theorem.
\end{abstract}

%
\begin{keyword}
\kwd{dislocations of wavefronts}
\kwd{random seas}
\kwd{Rice formulae}
\kwd{specular points}
\end{keyword}

\end{frontmatter}

\section{Introduction}

Many problems in applied mathematics require estimations of the number
of points, the length, the volume and so on, of the level sets of a
random function $\{ W(\mathbf{x})\dvtx\mathbf{x}\in\mathbb{R}^d \} $,
or of
some functionals defined on them. Let us mention some examples which
illustrate this general situation:
\begin{longlist}[1.]
\item[1.] A first example in dimension one is the number
of times that a random process
\\
$\{X(t)\dvtx t\in\mathbb{R}\}$ crosses the level
$u$:
\[
N_A^{X}(u)=\# \{ s \in A \dvtx X(s)=u\}.
\]
Generally speaking, the probability distribution of the random
variable $ N_A^{X}(u) $ is unknown, even for simple models of
the underlying process. However, there exist some formulae to
compute $\mathbb{E}( N_A^{X})$ and also higher order moments; see, for
example, \cite{crle}.

\item[2.] A particular case is the number of specular
points of a random curve or a random surface. Consider first the case
of a random curve. A light source placed at $(0,h_1)$
emits a ray that is reflected at the point $(x,W(x))$ of the curve
and the reflected ray is registered by an observer placed at
$(0,h_2)$. Using the equality between the angles of incidence and reflection
with respect to the normal vector to the curve (i.e.,
$N(x)=(-W'(x),1)$), an elementary computation gives
%
\begin{equation}\label{sp1}
W'(x) = \frac{\alpha_2 r_1 - \alpha_1 r_2}{x(r_2-r_1)},
\end{equation}
where $ \alpha_i := h_i -W(x)$ and $r_i := \sqrt{ x^2 +\alpha_i^2 }
$, $i=1,2$. The points $(x,W(x))$ of the curve such that $x$ is a
solution of
(\ref{sp1}) are called ``specular points''. For each Borel
subset $A$ of the real line, we denote by $\mathit{SP}_1(A)$
the number of specular points belonging to $A$. One of our aims is to
study the probability distribution of $\mathit{SP}_1(A)$.

\item[3.] The following approximation, which turns out to be very accurate
in practice for ocean waves, was introduced some time ago by
Longuet-Higgins (\cite{lo1lo1,lo2lo2}; see also \cite
{krakra}). If we suppose that $h_1$ and $h_2$ are large with respect
to $ W(x)$ and
$x$, then $ r_i = \alpha_i + x^2/( 2\alpha_i) + \mathrm{O}( h_i^{-3} )$.
(\ref{sp1}) can then be approximated by
%
\begin{equation}\label{sp2}
W'(x) \simeq\frac{ x}{2 } \frac{\alpha_1 + \alpha_2}{\alpha_1
\alpha_2}
\simeq\frac{ x}{2 }\frac{h_1 + h_2}{h_1 h_2} = k x,
\qquad \mbox{where }
k:=\frac{ 1}{2} \biggl( \frac{1}{h_1} +\frac{1}{h_2} \biggr).
\end{equation}
Set $Y(x) := W'(x) -kx$ and let $\mathit{SP}_2(A)$ denote the number of roots of
$Y(x)$ belonging to the set $A$, an approximation of $\mathit{SP}_1(A)$ under
this asymptotic. The first part of Section \ref{secspec} below will
be devoted to obtaining some results on the distribution of the random
variable $\mathit{SP}_2(\mathbb{R})$.

\item[4.] Let $W\dvtx Q \subset\mathbb{R}^d \to\mathbb{R}^{d'}$ with $d>d'$
be a random field and define the level set
\[
\mathcal{C}_Q^W(\mathbf{u})=\{\mathbf{x}\in Q\dvtx
W(\mathbf{x})=\mathbf{u}\}.
\]
Under certain general conditions, this set is a $(d-d')$-dimensional
manifold, but, in any case, its $(d-d')$-dimensional Hausdorff measure
is well defined. We denote this measure by $\sigma_{d-d'}$. Our
interest will be in computing the mean of the $\sigma_{d-d'}$-measure
of this level set, that is,
$\mathbb{E}[\sigma_{d-d'}(\mathcal{C}_Q^W(\mathbf{u}))]$, as well
as its
higher moments. It will also be of interest to compute
\[
\mathbb{E}\biggl[ \int_{\mathcal{C}_Q^W(\mathbf{u})}Y(s)\,\mathrm{d}\sigma_{d-d'}(s)
\biggr],
\]
where $Y(s)$ is some random field defined on the level set. One can
find formulae of this type, as well as a certain number of
applications, in
\cite{caca,wsws} ($d'=1$),
\cite{aw}, Chapter 6, and \cite{azlo}.

\item[5.] Another set of interesting problems is related to phase
singularities of random wavefronts. These correspond to
lines of darkness in light propagation, or threads of silence in
sound propagation \cite{be1be1}. In a mathematical framework, they
can be defined
as the locations of points where the amplitudes of waves vanish. If we
represent a wave as
\[
W(\mathbf{x},t)=\xi(\mathbf{x},t)+\mathrm{i}\eta(\mathbf{x},t),\qquad  \mathbf
{x}\in\mathbb{R}^d,
\]
where
$\xi, \eta$ are independent homogenous Gaussian random fields, then the
dislocations are the intersections of the two random surfaces
$\xi(\mathbf{x},t)=0,  \eta(\mathbf{x},t)=0$. Here, we only
consider the case $d=2 $. At fixed time, say $t=0$, we will compute the
expectation of the random variable $\#\{\mathbf{x}\in S: \xi(\mathbf
{x},0)=\eta(\mathbf{x},0)=0\}.$
\end{longlist}

The aim of this paper is threefold: (a) to re-formulate some known
results in
a modern language; (b) to prove a certain number of new results, both
for the exact and approximate models,
especially variance computations in cases in which only first moments
have been
known until now, thus contributing to improve the statistical methods
derived from the probabilistic results; (c) in some cases, to prove a
central limit theorem.

Rice formulae are our basic tools. For statements and proofs, we refer
to the recent book \cite{aw}. On the other hand, we are not giving
full proofs since the required computations are quite long and
involved; one can find details and some other examples that we do not
treat here in~\cite{azlomw}. For numerical computations, we use
MATLAB programs which are available at the site
\url{http://www.math.univ-toulouse.fr/\textasciitilde azais/prog/programs.html}.

In what follows, $\lambda_d$ denotes the Lebesgue measure in $\mathbb{R}^d$,
$\sigma_{d'}(B)$ the
$d'$-dimensional Hausdorff measure of a Borel set $B$ and $M^{\mathrm{T}}$ the
transpose of a matrix $M$. $(\mathit{const})$ is a positive constant whose
value may change from one occurrence to another. $p_{\xi}(x) $ is the
density of the random variable or vector $\xi$ at the point $x$,
whenever it exists. If not otherwise stated, all random fields are
assumed to be Gaussian and centered.

\section{Specular points in dimension one}\label{secspec}

\subsection{Expectation of the number of specular points}

We first consider the Longuet-Higgins approximation (\ref{sp2}) of the
number of SP $(x,W(x))$, that is,
\[
\mathit{SP}_2(I) = \# \{x \in I \dvtx Y(x) = W'(x) -kx =0 \}.
\]
We assume that $\{W(x)\dvtx x\in\mathbb{R}\}$ has $
\mathcal{C}^2 $ paths and is stationary. The Rice formula for the
first moment (\cite{aw}, Theorem 3.2) then applies and gives
\begin{eqnarray} \label{f:spec2}
\mathbb{E}(\mathit{SP}_2(I)) &=& \int_I\mathbb{E}\bigl(|Y'(x)| |Y(x) =0\bigr) p_{Y(x)}
(0)\, \mathrm{d}x = \int_I
\mathbb{E}(|Y'(x)|)\frac{ 1}{\sqrt{\lambda_2}} \varphi\biggl(\frac
{kx}{\sqrt
{\lambda_2}}\biggr)\, \mathrm{d}x\nonumber
\\[-8pt]\\[-8pt]
&=& \int_I
G\bigl(-k,\sqrt{\lambda_4}\bigr)\frac{1}{\sqrt{\lambda_2}} \varphi\biggl(\frac
{kx}{\sqrt{\lambda_2}}\biggr)\, \mathrm{d}x,\nonumber
\end{eqnarray}
where $\lambda_2$ and $\lambda_4$ are the spectral moments of $W$ and
%
\begin{eqnarray}\label{G}
G(\mu, \sigma) := \mathbb{E}(|Z|) ,\qquad  Z \sim N(\mu,\sigma^2) =\mu
[2\Phi
(\mu/\sigma)-1]+2 \sigma
\varphi(\mu/\sigma),
\end{eqnarray}
where $ \varphi(\cdot)$ and $\Phi(\cdot)$ are respectively the density and
cumulative distribution functions of the standard Gaussian
distribution.

If we look at the total number of
specular points over the whole line, we get
%
\begin{equation} \label{f:spectot}
\mathbb{E}(\mathit{SP}_2( \mathbb{R})) =\frac{ G( k, \sqrt{\lambda_4})}{k}
\simeq
\sqrt{\frac{2\lambda_4}{\uppi}}\frac{1}{k}\biggl(
1+\frac{1}{2}\frac{k^2}{\lambda_4}+
\frac{1}{24}\frac{k^4}{\lambda_4^2}+\cdots\biggr),
\end{equation}
which is the result given in \cite{lo1lo1}, part II, formula (2.14),
page 846.
Note that this quantity is an increasing function of $\frac{\sqrt
{\lambda_4} }{k}$.

We now turn to the computation of the expectation of the number of
specular points $\mathit{SP}_1(I)$ defined by
(\ref{sp1}).
It is equal to the number of zeros of the process
$ \{Z (x) := W'(x) -m_1(x,W(x))\dvtx x\in\mathbb{R}\},$ where
\[
m_1(x,w)=\frac{x^2-(h_1-w)(h_2-w)+\sqrt
{[x^2+(h_1-w)^2][x^2+(h_2-w)^2]}}{x(h_1+h_2-2w)}.
\]
Assume that the process $\{W(x)\dvtx x\in\mathbb{R}\}$ is Gaussian,
centered and
stationary, with $\lambda_0=1$. The process $Z$ is not
Gaussian, so we use \cite{aw}, Theorem 3.4, to get
\begin{eqnarray}\label{spec1}
\mathbb{E}( \mathit{SP}_1([a,b]))&=&\int_a^b \mathrm{d}x \int_{-\infty}
^{+\infty}\mathbb{E}\bigl( |Z'(x)||Z(x)=0,W(x)=w \bigr)\nonumber
\\[-8pt]\\[-8pt]
&&{}\qquad\hspace*{29pt} \times \frac{1}{\sqrt{2 \uppi}} \mathrm{e}^{-{w^2/2}} \frac{1}{\sqrt{2\uppi
\lambda_2}}\mathrm{e}^{-{m_1^2(x,w)/(2 \lambda_2)}}\, \mathrm{d}w.\nonumber
\end{eqnarray}
For the conditional expectation in (\ref{spec1}), note that
\[
Z'(x)=W''(x)-\frac{\partial m_1}{\partial
x}(x,W(x))-\frac{\partial m_1}{\partial w}(x,W(x))W'(x)
\]
so that under the condition $\{Z(x)=0,  W(x)=w\}$, we get
\[
Z'(x)=W''(x)-K(x,w),\qquad \mbox{where } K(x,w)=\frac{\partial
m_1}{\partial x}(x,w)+\frac{\partial m_1}{\partial
w}(x,w)m_1(x,w).
\]
Once again, using Gaussian regression, we can write (\ref{spec1}) in the form
%
\begin{eqnarray}\label{spec3}
\mathbb{E}( \mathit{SP}_1([a,b]))
=\frac{1}{2
\uppi}\sqrt{\frac{\lambda_4-\lambda_2^2}{\lambda_2}}\int_a^b \mathrm{d}x
\int
_{-\infty} ^{+\infty}
G(m,1)
\exp\biggl(-\frac{1}{2}\biggl(w^2+\frac{m_1^2(x,w)}{\lambda
_2}\biggr)\biggr)\, \mathrm{d}w,\quad
\end{eqnarray}
where $ m=m(x,w)=(\lambda_2w+K(x,w))/\sqrt{\lambda_4-\lambda
_2^2}$ and $G$ is defined in (\ref{G}). In (\ref{spec3}), the
integral is convergent as
$a\rightarrow-\infty,b\rightarrow+\infty$ and this formula
is well adapted to numerical approximation.

\begin{figure}[t]

\includegraphics{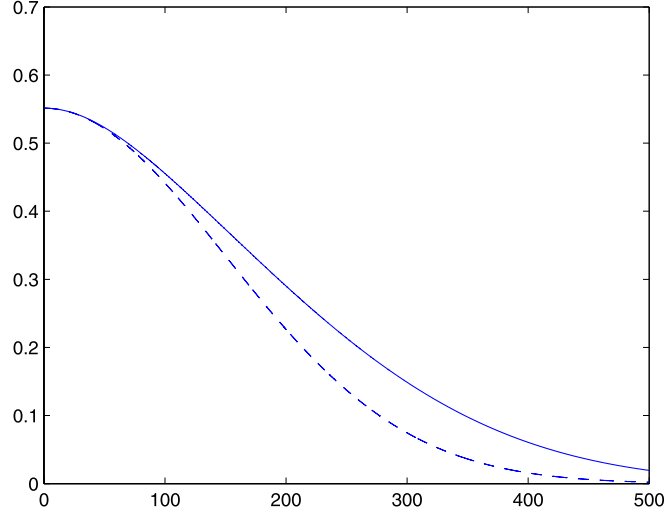}

\caption{ Intensity of specular points in the case $h_1= 100,h_2= 300,
\lambda_4=3$.
Solid line corresponds to the exact formula, dashed line corresponds to
the approximation (\protect\ref{f:spec2}).}\label{fi:1}
\end{figure}

We have performed some numerical computations to compare
the exact expectation given by (\ref{spec3}) with the
approximation (\ref{f:spec2}) in the stationary case. The result
depends on $ h_1,h_2, \lambda_4$ and
$\lambda_2$, and, after scaling, we can assume that $\lambda_2=1$.
When $h_1 \approx h_2$, the approximation
(\ref{f:spec2}) is very sharp. For example, if $h_1=100,h_2=100,
\lambda_4=3$,
the expectation of the total number of specular
points over $ \mathbb{R}$ is $138.2$; using the approximation (\ref
{f:spectot}),
the result with the exact formula is around $ 2 .10^{-2}$ larger (this
is the same order as the error in the
computation of the integral). For $h_1=90,h_2=110, \lambda_4=3$, the
results are $136.81$ and $137.7$, respectively.
If $h_1=100,h_2=300, \lambda_4=3$, the results differ significantly and
Figure \ref{fi:1} displays the densities in the integrand of (\ref
{spec1}) and
(\ref{f:spec2}) as functions of $x$.

\subsection{Variance of the number of specular points}

We assume that the covariance function $\mathbb{E}(W(x)W(y))=\Gamma
(x-y)$ has enough regularity to perform
the computations below, the precise requirements being given
in the statement of Theorem~\ref{asymptoticvar}.

Writing, for short, $S=\mathit{SP}_2( \mathbb{R})$, we have
%
\begin{equation}\label{vardesc}
\operatorname{Var}(S)=\mathbb{E}\bigl(S(S-1)\bigr)+\mathbb{E}(S)-[\mathbb{E}(S)]^2.
\end{equation}
Using \cite{aw}, Theorem 3.2, we have
\begin{eqnarray}\label{ricesp2}
\mathbb{E}\bigl(S(S-1)\bigr)&=&\int_{\mathbb{R}^2} \mathbb{E}
\bigl(|W''(x)-k||W''(y)-k||
W'(x)=kx,W'(y)=ky
\bigr) \nonumber
\\[-8pt]\\[-8pt]
&&{}\hspace*{12pt}\times p_{W'(x),W'(y)}(kx,ky)\,\mathrm{d}x\,\mathrm{d}y,\nonumber
\end{eqnarray}
where
\begin{eqnarray}\label{densite}
&&p_{W'(x),W'(y)}(kx,ky)\nonumber
\\[-8pt]\\[-8pt]
&&\quad =\frac{1}{2\uppi
\sqrt{\lambda_2^2-\Gamma''^2(x-y)}}\exp\biggl[-\frac{1}{2}
\frac{k^2(\lambda_2x^2+2\Gamma''^2(x-y)xy+\lambda_2y^2)}{\lambda
_2^2-\Gamma''^2(x-y)}
\biggr],\nonumber
\end{eqnarray}
under the condition that the
density (\ref{densite}) does not degenerate for $x\neq y$.

For the
conditional expectation in (\ref{ricesp2}), we perform a Gaussian
regression of $W''(x)$ (resp., $W''(y)$) on the pair $(W'(x),W'(y))$.
Putting $z=x-y$, we obtain
\begin{eqnarray*}
W''(x)&=&\theta_y(x)+a_y(x)W'(x)+b_y(x)W'(y),
\\
a_y(x)&=&-\frac{\Gamma'''(z)\Gamma''(z)}{\lambda_2^2-\Gamma
''^2(z)}, \qquad b_y(x)=-\frac{\lambda_2\Gamma'''(z)}{\lambda
_2^2-\Gamma''^2(z)},
\end{eqnarray*}
where $\theta_y(x)$ is Gaussian centered, independent of
$(W'(x),W'(y))$. The regression of $W''(y) $ is obtained by
permuting $x$ and $y$.

The conditional expectation in (\ref{ricesp2}) can now be rewritten
as an unconditional expectation:
%
\begin{equation}\label{cond1}
\mathbb{E}\biggl\{ \biggl|\theta_y(x)-k\Gamma'''(z)\biggl[
1+\frac{\Gamma''(z)x+\lambda_2y}{\lambda_2^2-\Gamma''^2(z)}\biggr]
\biggr| \biggl|\theta_x(y)-k\Gamma'''(-z)\biggl[
1+\frac{\Gamma''(-z)y+\lambda_2x}{\lambda_2^2-\Gamma''^2(z)}\biggr]
\biggr|\biggr\}.\
\end{equation}
Note that the singularity on the diagonal $x=y$ is removable
since a Taylor expansion shows that for $z\approx0$,
%
\begin{equation}\label{eq1}
\Gamma^{\prime\prime\prime} (z)\biggl[
1+\frac{\Gamma''(z)x+\lambda_2y}{\lambda_2^2-\Gamma''^2(z)}\biggr]=
\frac{1}{2}\frac{\lambda_4}{\lambda_2}x \bigl(z+\mathrm{O}(z^3) \bigr).
\end{equation}
It can be checked that
%
\begin{eqnarray}\label{cond2}
\sigma^2(z)&=&\mathbb{E}( (\theta_y(x))^2 )=\mathbb{E}(
(\theta_x(y))^2
)=\lambda_4-\frac{\lambda_2\Gamma'''^2(z)}{\lambda_2^2-\Gamma
''^2(z)},
\\\label{cond3}
\mathbb{E}( \theta_y(x)\theta_x(y)
)&=&\Gamma^{(4)}(z)+\frac{\Gamma'''^2(z)\Gamma''(z)}{\lambda
_2^2-\Gamma''^2(z)}.
\end{eqnarray}
Moreover, if $\lambda_6<+\infty$,
we can show that as $z\approx0$, we have
%
\begin{equation}\label{eq2}
\sigma^2(z)\approx\frac{1}{4}\frac{\lambda_2\lambda_6-\lambda
_4^2}{\lambda_2}z^2
\end{equation}
and it follows that the singularity on the diagonal of the
integrand
in the right-hand side of (\ref{ricesp2}) is also removable.

We will make use of the following auxiliary statement that we state
as a lemma for further reference. The proof requires some
calculations, but is elementary, so we omit it. The value of
$H(\rho;0,0)$ can be found in, for example, \cite{crle}, pages
211--212.
\begin{lem}\label{esproabsv}
Let
\[
H(\rho;\mu,\nu)=\mathbb{E}(|\xi+\mu||\eta+\nu|),
\]
where the pair $(\xi,\eta)$ is centered Gaussian, $\mathbb{E}(\xi
^2)=\mathbb{E}(\eta^2)=1,\mathbb{E}(\xi\eta)=\rho.$

Then, if $\mu^2+\nu^2 \leq1$ and $0\leq\rho\leq1$,
\[
H(\rho;\mu,\nu)= H(\rho;0,0)+R_2(\rho;\mu,\nu),
\]
where
\[
H(\rho;0,0)=\frac{2}{\uppi}\sqrt{1-\rho^2}+\frac{2\rho}{\uppi
}\arctan
\frac{\rho}{\sqrt{1-\rho^2}}\quad \mbox{and}\quad |R_2(\rho;\mu,\nu
)|\leq3 (\mu^2+\nu^2).
\]
\end{lem}

In the next theorem, we compute the equivalent of the variance of the
number of specular points, under certain hypotheses on the random
process $W$ and with the Longuet-Higgins asymptotic. This result is new
and useful for estimation purposes since it implies that, as
$k\rightarrow0$, the coefficient of variation of the random
variable $S$ tends to zero at a known speed. Moreover, it will also
appear in a natural way when normalizing $S$ to obtain a central
limit theorem.

\begin{theo}\label{asymptoticvar}
Assume that the centered Gaussian stationary process
$\mathcal{W}=\{W(x) \dvtx x\in\mathbb{R}\}$ is $\delta$-dependent, that is,
$\Gamma(z)=0$ if $|z|>\delta$, and that it has
$\mathcal{C}^4$-paths. Then, as $k\rightarrow0$, we have
%
\begin{equation}\label{eqvar}
\operatorname{Var}(S)=\theta\frac{1}{k}+\mathrm{O}(1),
\end{equation}
where
\begin{eqnarray*}
\theta&=&\Biggl(\frac{J}{\sqrt{2}}+\sqrt{\frac{2\lambda
_4}{\uppi}}-\frac{2\delta\lambda_4}{\sqrt{\uppi^3 \lambda_2}} \Biggr)
, \qquad
J =\int_{-\delta}^{+\delta}\frac{\sigma^2(z)H(\rho(z);0,0)
)}{\sqrt{2\uppi(\lambda_2+\Gamma''(z))}}\,\mathrm{d}z,
\\
\rho(z) &=&\frac{1}{\sigma^2(z)}\biggl[ \Gamma^{(4)}(z)+\frac{\Gamma
'''(z)^2\Gamma''(z)}{\lambda_2^2-\Gamma''^2(z)}\biggr],
\end{eqnarray*}
$\sigma^2(z)$ is defined in (\ref{cond2}) and
$H$ is defined in Lemma \ref{esproabsv}.
Moreover, as
$k\rightarrow0$, we have
\[
\frac{\sqrt{\operatorname{Var}(S)}}{\mathbb{E}(S)}\approx\sqrt
{\theta k}.
\]
\end{theo}

\begin{Rem*}
\begin{enumerate}[(1)]
\item[(1)] The $\delta$-dependence hypothesis can be replaced
by some weaker mixing condition, such as
\[
\bigl|\Gamma^{(i)}(z) \bigr| \leq (\mathit{const})(1+|z|)^{-\alpha
}\qquad (0\leq
i \leq4)
\]
for some $\alpha>1$, in which case the value of $\theta$ should be
\[
\theta=\sqrt{\frac{2\lambda_4}{\uppi}}+\frac{1}{\sqrt{\uppi}}\int
_{-\infty}^{+\infty}
\biggl[ \frac{\sigma^2(z)H(\rho(z);0,0)}{2\sqrt{\lambda_2+\Gamma
''(z)}}-\frac{1}{\uppi}\frac{\lambda_4}{\sqrt{\lambda_2}}\biggr]\,\mathrm{d}z.
\]
The proof of this extension can be constructed along the same
lines as the one we give below, with some additional computations.

\item[(2)] The above computations complete the study done  in \cite
{lo1lo1} (Theorem 4).
In \cite{krakra}, the random variable
$\mathit{SP}_2(I)$ is expanded in the Wiener--Hermite chaos. The aforementioned
expansion yields
the same formula for the expectation and also allows a
formula to be obtained for the variance. However, this expansion is
difficult to manipulate in
order to get the result of Theorem \ref{asymptoticvar}.
\end{enumerate}
\end{Rem*}

\begin{pf*}{Proof of Theorem \ref{asymptoticvar}}
We use the notation and the computations preceding the statement of
the
theorem.

Divide the integral on the right-hand side of (\ref{ricesp2}) into
two parts, corresponding to $|x-y|>\delta$ and $|x-y|\leq\delta$,
that is,
%
\begin{equation}\label{dividiren2}
\mathbb{E}\bigl(S(S-1)\bigr)=\int\hspace*{-2pt}\int_{|x-y|>\delta}\cdots+\int\hspace*{-2pt}\int
_{|x-y|\leq
\delta}\cdots=I_1+I_2.
\end{equation}
In the first term, the $\delta$-dependence of the process implies
that one can factorize the conditional expectation and the density
in the integrand. Taking into account that for each $x\in\mathbb{R}$, the
random variables $W''(x)$ and $W'(x)$ are independent, we obtain for
$I_1$\vspace*{2pt}
\begin{eqnarray*}
I_1=\int\hspace*{-2pt}\int_{|x-y|>\delta}\mathbb{E}\bigl(|W''(x)-k| \bigr)\mathbb{E}
\bigl(|W''(y)-k| \bigr)p_{W'(x)}(kx)p_{W'(y)}(ky)\,\mathrm{d}x\,\mathrm{d}y.\vspace*{2pt}
\end{eqnarray*}
On the other hand, we know that $W'(x)$ (resp., $W''(x)$) is
centered normal with variance $\lambda_2$ (resp., $\lambda_4$).
Hence,\vspace*{2pt}
\[
I_1=\bigl[G\bigl(k,\sqrt{\lambda_4}\bigr)
\bigr]^2\int\hspace*{-2pt}\int_{|x-y|>\delta}\frac{1}{2\uppi\lambda_2}\exp\biggl[
-\frac{1}{2}\frac{k^2(x^2+y^2)}{\lambda_2} \biggr]\,\mathrm{d}x\,\mathrm{d}y.\vspace*{2pt}
\]
To compute the integral on the right-hand side, note that the
integral over the whole $x,y$ plane is equal to $1/k^2$ so that it
suffices to compute the integral over the set ${|x-y|\leq\delta}$.
Changing variables, this last integral is equal to\vspace*{2pt}
\begin{eqnarray*}
\int_{-\infty}^{+\infty}\mathrm{d}x\int_{x-\delta}^{x+\delta}\frac
{1}{2\uppi
\lambda_2}\exp\biggl[ -\frac{1}{2}\frac{k^2(x^2+y^2)}{\lambda_2}
\biggr]\,\mathrm{d}y =\frac{\delta}{k\sqrt{\lambda_2\uppi}}+\mathrm{O}(1),\vspace*{2pt}
\end{eqnarray*}
where the last term is bounded if $k$ is bounded (remember
that we are considering an approximation in which $k\approx0$). Therefore,
we can conclude that
\[
\int\hspace*{-2pt}\int_{|x-y|>\delta}\frac{1}{2\uppi\lambda_2}\exp\biggl[
-\frac{1}{2}\frac{k^2(x^2+y^2)}{\lambda_2}
\biggr]\,\mathrm{d}x\,\mathrm{d}y=\frac{1}{k^2}-\frac{\delta}{k\sqrt{\lambda_2\uppi}}+\mathrm{O}(1),
\]
from which we deduce, performing a Taylor
expansion, that
%
\begin{equation}\label{i1}
I_1=\frac{2\lambda_4}{\uppi}\biggl[
\frac{1}{k^2}-\frac{\delta}{k\sqrt{\lambda_2 \uppi}}+\mathrm{O}(1)\biggr].
\end{equation}
Let us now turn to $I_2$. Using Lemma \ref{esproabsv} and the
equivalences (\ref{eq1}) and
(\ref{eq2}), whenever $|z|=|x-y|\leq\delta$, the integrand on the
right-hand side of (\ref{ricesp2}) is bounded by
\[
(\mathit{const})[ H(\rho(z);0,0)+k^2(x^2+y^2)].
\]
We divide the integral $I_2$ into two parts.

First, on the set $\{(x,y)\dvtx |x|\leq2\delta, |x-y|\leq
\delta\}$, the integral is clearly bounded by some constant.

Second, we consider the integral on the set $\{(x,y)\dvtx x> 2\delta,
|x-y|\leq\delta\}$. (The symmetric case, replacing $x> 2\delta$ by
$x<- 2\delta$, is similar -- that is the reason for the factor $2$ in
what follows.) We have (recall that $z=x-y$)
\begin{eqnarray*}
I_2&=&\mathrm{O}(1)+2\int\hspace*{-2pt}\int_{|x-y|\leq\delta,x>2\delta}\sigma^2(z)[
H(\rho(z);0,0 )+R_2(\rho(z);\mu,\nu)]
\\
&&{}\qquad\hspace*{79pt} \times\frac{1}{2\uppi\sqrt{\lambda_2^2-\Gamma''^2(z)}}
\\
&&{}\qquad\hspace*{79pt}\times \exp
\biggl[-\frac{1}{2}\frac{k^2(\lambda_2x^2+2\Gamma''(x-y)xy+\lambda
_2y^2)}{\lambda_2^2-\Gamma''^2(x-y)}\biggr]\,\mathrm{d}x\,\mathrm{d}y,
\end{eqnarray*}
which can be rewritten as
\begin{eqnarray*}
I_2&=&\mathrm{O}(1)+2\int_{-\delta} ^{\delta}\sigma^2(z)[
H(\rho(z);0,0 )+R_2(\rho(z);\mu,\nu)]
\\
&&{}\qquad\hspace*{31pt} \times\frac{1}{\sqrt{2\uppi(\lambda_2+\Gamma''(z))}}\exp
\biggl[-\frac{1}{2}\frac{k^2z^2}{\lambda_2-\Gamma''(z)}
\biggl(\frac{\lambda_2}{\lambda_2+\Gamma''(z)}-\frac{1}{2}
\biggr)\biggr]\,\mathrm{d}z
\\
&&{}\qquad\hspace*{31pt} \times\int_{2\delta}^{+\infty}\frac{1}{\sqrt{2\uppi(\lambda
_2-\Gamma''(z))}}\exp\biggl[-k^2\frac{(x-z/2)^2}{\lambda_2-\Gamma
''(z)} \biggr]\,\mathrm{d}x.
\end{eqnarray*}
Changing variables, the inner integral becomes
%
\begin{equation}\label{intint}
\frac{1}{k\sqrt{2}}\int_{\tau_0}^{+\infty}\frac{1}{\sqrt{2\uppi
}}\exp\biggl( -\frac{1}{2}\tau^2
\biggr)\,\mathrm{d}\tau=\frac{1}{2\sqrt{2}}\frac{1}{k}+\mathrm{O}(1),
\end{equation}
where
$\tau_0=2\sqrt{2}k(2\delta-z/2)/\sqrt{\lambda_2-\Gamma''(z)}$.

Substituting this into $I_2$, we obtain
%
\begin{equation}\label{finali2}
I_2=\mathrm{O}(1)+\frac{J}{k\sqrt{2}}.
\end{equation}
To finish, combine (\ref{finali2}) with (\ref{i1}),
(\ref{dividiren2}), (\ref{vardesc}) and (\ref{f:spectot}).
\end{pf*}

\subsection{Central limit theorem}
\begin{theo}\label{tcl}
Assume that the process $\mathcal{W} $ satisfies the hypotheses of
Theorem \ref{asymptoticvar}. In addition, we assume that the
fourth moment of the number of approximate specular points on an
interval having length equal to $1$ is uniformly bounded in $k$,
that is, for all $a\in\mathbb{R}$ and $0<k<1$,
%
\begin{equation}\label{mom4}
\mathbb{E}\bigl(\bigl[\mathit{SP}_2([a,a+1]) \bigr]^4\bigr)\leq(\mathit{const}).
\end{equation}
Then, as $ k\rightarrow0$,
\begin{eqnarray*}
\frac{S - \sqrt{{2\lambda_4/\uppi}}{1/k}}{\sqrt{\theta
/k}} \quad\Longrightarrow\quad N(0,1)\qquad \mbox{in distribution.}
\end{eqnarray*}
\end{theo}

\begin{Rem*} One can give conditions for the additional
hypothesis (\ref{mom4}) to
hold true. Even though they are not nice, they are not costly from the
point of
view of physical models. For example, either one of the following
conditions implies (\ref{mom4}):
\begin{longlist}[(ii)]
\item[(i)] the paths $x\rightsquigarrow W(x)$ are of class $\mathcal
{C}^{11}$ (use \cite{aw}, Theorem 3.6,
with $m=4$, applied to the random process $\{W'(x)\dvtx x\in\mathbb{R}\}
$);

\item[(ii)] the paths $x\rightsquigarrow W(x)$ are of class $\mathcal
{C}^{9}$ and the support of the spectral
measure has an accumulation point (apply \cite{aw}, Example 3.4,
Proposition 5.10 and Theorem 3.4, to show
that the fourth moment of the number of zeros of $W''(x)$ is
bounded).
\end{longlist}
Note that the asymptotic here differs from other ones existing in the
literature on related subjects (compare with, e.g., \cite{cuzcuz} and
\cite{pirypiry}).
\end{Rem*}

\begin{pf*}{Proof of Theorem \ref{tcl}} Let $\alpha$ and $\beta$ be real numbers
satisfying the conditions
$1/2<\alpha<1$, $\alpha+ \beta>1$, $2\alpha+ \beta<2$.
It suffices to prove the convergence as $k$ takes values on a
sequence of positive numbers tending to $0$. To keep in mind that the
parameter is $k$,
we use the notation $S(k):=S=\mathit{SP}_2(\mathbb{R}) $.

Choose $k$ small enough so that $k^{-\alpha} >2$ and define
the sets of disjoint intervals, for $j=0, \pm1,\ldots, \pm[k^{-\beta
}]$ ([$\cdot$] denotes integer part),\vspace*{-2pt}
\begin{eqnarray*}
U_j^k &=& \bigl((j-1) [k ^{-\alpha}]\delta+\delta/2, j [k ^{-\alpha
}]\delta-\delta/2\bigr),
\\[-2pt]
I_j^k&=&\bigl[j[k ^{-\alpha}]\delta-\delta/2,j[k
^{-\alpha}]\delta+\delta/2 \bigr].\vspace*{-2pt}
\end{eqnarray*}
Each interval $U_j^k$ has length $[k ^{-\alpha}]\delta
-\delta$ and two neighboring intervals $U_j^k$ are separated by
an interval of length $\delta$. So, the $\delta$-dependence of the
process implies that the random variables $\mathit{SP}_2(U_j^k),j=0, \pm
1,\ldots, \pm[k^{-\beta}]$, are independent. A similar argument
applies to $\mathit{SP}_2(I_j^k),j=0, \pm
1,\ldots, \pm[k^{-\beta}]$.

We write\vspace*{-2pt}
\[
T(k) = \sum_{|j|\leq[k ^{-\beta}]} \mathit{SP}_2( U_j^k),\qquad V_k =
(\operatorname{Var}
( S(k)) ) ^{-1/2} \approx\sqrt{k/\theta},\vspace*{-2pt}
\]
where the equivalence is due to Theorem \ref{asymptoticvar}.

The proof is performed in two steps, which easily imply the statement. In
the first, it is proved that $ V_k [ S(k) -T(k)]$ tends to $0$ in the
$L^2$ of the underlying probability space. In the second step, we prove
that $ V_kT(k)$ is asymptotically standard normal.

\textit{Step} 1. We first prove that $V_k [S(k) -T(k)] $ tends to
$0$ in $L^1$. Since it is non-negative, it suffices to show that its
expectation tends to zero. We have
\[
S(k) -T(k)= \sum_{|j|< [k ^{-\beta}]}\mathit{SP}_2 (I_j^k)+Z_1+Z_2,
\]
where
$Z_1=\mathit{SP}_2 (-\infty,- [k ^{-\beta}]\cdot [k ^{-\alpha}]\delta
+\delta
/2 )$,
$Z_2=\mathit{SP}_2 ( [k ^{-\beta}]\cdot [k ^{-\alpha}]\delta-\delta
/2,+\infty)).$

 Using the fact that $\mathbb{E}(\mathit{SP}^k_2(I)) \leq(const)
\int_I
\varphi( kx/\sqrt{\lambda_2}) \,\mathrm{d}x$,  we can show that
\[
V_k \mathbb{E}\bigl(S(k) -T(k)\bigr) \leq(\mathit{const})k^{1/2} \Biggl[ \sum_{\ell
=0}^{+\infty} \varphi\biggl(\frac{\ell[k^{-\alpha}]k\delta}{ \sqrt
{\lambda_2}}\biggr)
+ \int_{ [k^{-\alpha}][k^{-\beta}]\delta} ^{+\infty} \varphi\bigl(
kx/\sqrt{\lambda_2}\bigr) \,\mathrm{d}x \Biggr],
\]
which tends to zero as a consequence of the choice of $\alpha$ and
$\beta$.
It suffices to prove that $V_k^2 \operatorname{Var}( S(k) -T(k)
)
\rightarrow0 $ as $k\rightarrow0$. Using independence, we have
\begin{eqnarray*}
\operatorname{Var}\bigl( S(k) -T(k) \bigr)&=& \sum_{|j|< [k ^{-\beta
}]}\operatorname{Var}(
\mathit{SP}_2 (I_j^k) ) + \operatorname{Var}(Z_1)+ \operatorname{Var}(Z_2
)\\
&\leq&\sum_{|j|< [k ^{-\beta}]}\mathbb{E}\bigl( \mathit{SP}_2 (I_j^k)\bigl(\mathit{SP}_2
(I_j^k)-1\bigr) \bigr)
\\
&&{}+ \mathbb{E}\bigl(Z_1(Z_1-1)\bigr)+ \mathbb{E}\bigl(Z_2(Z_2-1) \bigr)+\mathbb
{E}\bigl( S(k)-T(k) \bigr).
\end{eqnarray*}
We already know that $V_k^2\mathbb{E}( S(k)-T(k)
)\rightarrow
0.$ Since each $I_j^k$ can be covered by a fixed\vspace*{-2pt} number of intervals of
size one, we know that
$\mathbb{E}( \mathit{SP}_2 (I_j^k)(\mathit{SP}_2 (I_j^k)-1) ) $ is
bounded by a constant which does not depend on $k$ and $j$. Therefore,
\[
V_k^2\sum_{|j|< [k ^{-\beta}]}\mathbb{E}\bigl( \mathit{SP}_2
(I_j^k)\bigl(\mathit{SP}_2 (I_j^k)-1\bigr) \bigr)\leq(\mathit{const})k^{1-\beta},
\]
which tends to zero because of the choice of $\beta$.
The remaining two terms can be bounded in a similar form as in the
proof of Theorem \ref{asymptoticvar}.

\textit{Step} 2. $T(k)$ is a sum of independent,
but not equidistributed, random variables. To prove that it satisfies a
central limit theorem,
we will use a Lyapunov condition based of fourth moments. Set
\[
M_j^m := \mathbb{E}\{[ \mathit{SP}_2(U_j^k) -\mathbb{E}(
\mathit{SP}_2(U_j^k))
]^m\}.
\]
For the Lyapunov condition, it suffices to verify that
%
\begin{equation}\label{lyap}
\Sigma^{-4}\sum_{|j|\leq[k ^{-\beta}]} M_j^4 \to0\qquad \mbox{as }
k\rightarrow0,
\mbox{ where }
\Sigma^2 := \sum_{|j|\leq[k ^{-\beta}]} M_j^2.
\end{equation}
To prove (\ref{lyap}), we divide each interval $U_j^k$ into
$ p=[k^{-\alpha}]-1$ intervals $I_1,\dots,I_p$ of equal size~$\delta$.
We have
%
\begin{equation}\label{sumamom4}
\mathbb{E}(\mathit{SP}_1+\cdots+\mathit{SP}_p)^4 = \sum_{1\leq i_1,i_2,i_3,i_4
\leq p} \mathbb{E}
( \mathit{SP}_{i_1}\mathit{SP}_{i_2} \mathit{SP}_{i_3} \mathit{SP}_{i_4}),
\end{equation}
where $ \mathit{SP}_i$ stands for $\mathit{SP}_2(I_i) -\mathbb{E}( \mathit{SP}_2(I_i)
)$. Since
the size of all intervals
is equal to $\delta$, given the finiteness of fourth moments in the
hypothesis, it follows that
$\mathbb{E}( \mathit{SP}_{i_1}\mathit{SP}_{i_2} \mathit{SP}_{i_3} \mathit{SP}_{i_4})$ is
bounded.

On the other hand, the number of terms which do not
vanish in the sum of the right-hand side of (\ref{sumamom4}) is
$\mathcal{ O}(p^2 )$. In fact, if one of the indices in
$(i_1,i_2,i_3,i_4)$ differs by more than $1$ from all the others, then
$\mathbb{E}( \mathit{SP}_{i_1}\mathit{SP}_{i_2} \mathit{SP}_{i_3} \mathit{SP}_{i_4})=0$. Hence,
\[
\mathbb{E}[ \mathit{SP}_2(U_j^k) -\mathbb{E}( \mathit{SP}_2(U_j^k))
]^4 \leq(const)
k^{-2\alpha}
\]
so that $\sum_{|j|\leq[k^{-\beta}]} M_j^4 = \mathcal{O}
(k^{-2\alpha} k^{-\beta}).$ The inequality $2 \alpha+\beta<2 $
implies the Lyapunov condition.
\end{pf*}

\section{Specular points in two dimensions. Longuet-Higgins
approximation} \label{s:spec2d}

We consider, at fixed time, a random surface depending on two space
variables $x$ and $y$. The source of light is placed at $(0,0,h_1)$
and the observer is at $(0,0, h_2)$. The point $(x,y)$ is a specular
point if the normal vector $ n(x,y) = ( -W_x,-W_y,1)$ to the surface
at $(x,y)$ satisfies the following two conditions:
\begin{itemize}[$\bullet$]
\item the angles with the incident ray $I = (
-x,-y,h_1-W)$ and the reflected ray $R = (
-x,-y,h_2-W)$ are equal (to simplify notation, the argument $(x,y)$ has
been removed);
\item it belongs to the plane generated by $I$ and $R$.
\end{itemize}
Setting $ \alpha_i =h_i -W$ and $ r_i = \sqrt{x^2 +y^2 +\alpha_i}$,
$i=1,2$, as in the one-parameter case, we have
%
\begin{equation}\label{spec:2:1}
W_x = \frac{x}{x^2 +y^2}\frac{ \alpha_2 r_1 - \alpha_1
r_2}{r_2-r_1},\qquad
W_y = \frac{y}{x^2 +y^2}\frac{ \alpha_2 r_1 - \alpha_1 r_2}{r_2-r_1}.
\end{equation}
When $h_1$ and $h_2$ are large, the system above can be approximated by
%
\begin{equation}\label{spec:2:2}
W_x = kx,\qquad  W_y = ky,
\end{equation}
under the same conditions as in dimension one.

Next, we compute the expectation of $ \mathit{SP}_2(Q)$, the number of
approximate specular points, in the sense of (\ref{spec:2:2}), that
are in a domain $Q$. In the remainder of this paragraph, we limit our
attention to this approximation and to the case in which $\{
W(x,y)\dvtx (x,y)\in\mathbb{R}^2 \}$ is a centered
Gaussian stationary random field.

Let us define
%
\begin{equation}\label{f:specd}
\mathbf{Y}(x,y) :=
\pmatrix{
W_x(x,y) -kx \cr
W_{y} (x,y)-ky
}
.
\end{equation}
Under very general conditions, for example, on the spectral measure
of $\{W(x,y)\dvtx x,y\in\mathbb{R}\}$, the random field $\{Y(x,y)\dvtx x,y\in
\mathbb{R}\}$
satisfies the conditions of \cite{aw}, Theorem 6.2, and we can write
%
\begin{equation}\label{f:spec22}
\mathbb{E}( \mathit{SP}_2(Q)) = \int_Q \mathbb{E}( |\det
\mathbf
{Y}'(x,y)|)
p_{\mathbf{Y}(x,y)}(\mathbf{0})\, \mathrm{d}x \,\mathrm{d}y
\end{equation}
since for fixed $(x,y)$, the random matrix $\mathbf{Y}'(x,y)$ and the random
vector $\mathbf{Y}(x,y)$ are independent so that the condition in the
conditional expectation can be eliminated. The density in the
right-hand side of (\ref{f:spec22}) has the expression
\begin{eqnarray}\label{f:d0}
p_{\mathbf{Y}(x,y)}(\mathbf{0}) &=& p_{ (W_x,W_y)} (kx,ky)\nonumber
\\[-8pt]\\[-8pt]
&=&
\frac{1}{2\uppi}\frac{1}{\sqrt{\lambda_{20}\lambda_{02}-\lambda_{11}^2}}
 \exp\biggl[
-\frac{k^2}{2(\lambda_{20}\lambda_{02}-\lambda_{11}^2)}
(\lambda_{02}
x^2-2\lambda_{11}xy+\lambda_{20}y^2 )\biggr].\nonumber\qquad
\end{eqnarray}
To compute the expectation of the absolute value of the determinant
in the right-hand side of (\ref{f:spec22}), which does not depend on
$x,y$, we
use the method of \cite{be1be1}. Set $ \Delta:=
\det\mathbf{Y}'(x,y)=( W_{xx} -k) (W_{yy} -k) -
W^2_{xy}$.

We have
%
\begin{equation} \label{f:d1}
\mathbb{E}(|\Delta|) = \mathbb{E}\biggl[\frac{2}{\uppi} \int_0
^{+\infty} \frac{1-
\cos( \Delta t)}{t^2}\, \mathrm{d}t \biggr].
\end{equation}
Define
\[
h(t) := \mathbb{E}\bigl[ \exp\bigl( \mathrm{i}t [ ( W_{xx} -k) ( W_{yy} -k)
-W_{xy}^2] \bigr)\bigr].
\]
Then
%
\begin{equation}\label{f:d2}
\mathbb{E}(|\Delta|) =\frac{2}{\uppi} \biggl( \int_0 ^{+\infty}
\frac{1-
\mathfrak{Re} [h(t)]}{t^2} \,\mathrm{d}t \biggr).
\end{equation}
We now proceed to give a formula for $\mathfrak{Re} [h(t)]$. Define
\[
A =
\pmatrix{
0 & 1/2& 0\cr
1/2 & 0 & 0\cr
0& 0 & -1
}
\]
and denote by $ \Sigma$ the variance matrix of $ ( W_{xx},
W_{yy},W_{x,y} )$
\[
\Sigma: =
\pmatrix{
\lambda_{40} & \lambda_{22}& \lambda_{31}\cr
\lambda_{22} & \lambda_{04} & \lambda_{13}\cr
\lambda_{31}& \lambda_{13} & \lambda_{22}
}.
\]
Let $\Sigma^{1/2} A\Sigma^{1/2} = P
\operatorname{diag}(\Delta_1,\Delta_2,\Delta_3) P^{\mathrm{T}}$, where $P$ is orthogonal. Then
\begin{eqnarray} \label{m:1}
h(t) &=& \mathrm{e}^{\mathrm{i}tk^2}
\mathbb{E}\bigl( \exp\bigl[ \mathrm{i}t \bigl( \bigl( \Delta_1 Z^2_1 -k (s_{11} + s_{21})
Z_1 \bigr)
+ \bigl( \Delta_2 Z^2_2 -k (s_{12} + s_{22}) Z_2 \bigr)\nonumber
\\[-8pt]\\[-8pt]
&&{}\hspace*{56pt}+\bigl( \Delta_3 Z^2_3 -k (s_{13} + s_{23}) Z_3 \bigr)\bigr) \bigr] \bigr),\nonumber
\end{eqnarray}
where $ (Z_1,Z_2,Z_3)$ is standard normal and $ s_{ij} $ are the
entries of
$\Sigma^{1/2} P^\mathrm{T}$.

One can check that if $ \xi$ is a standard normal variable and
$\tau, \mu$ are real constants, $\tau>0$, then
\[
\mathbb{E}\bigl( \mathrm{e}^{ \mathrm{i}\tau( \xi+ \mu)^2 } \bigr) =(1-2\mathrm{i}\tau
)^{-1/2}\mathrm{e}^{{\mathrm{i}\tau\mu^2/(1-2\mathrm{i}\tau)}}
= \frac{1}{(1+4\tau^2)^{1/4}}
\exp\biggl[\frac{-2\tau}{1+4\tau^2} +\mathrm{i}\biggl(\varphi+\frac{\tau\mu
^2}{1+4\tau^2}
\biggr)\biggr],
\]
where $ \varphi=\frac{1}{2}\arctan(2\tau),0<\varphi<\uppi/4.$
Substituting this into (\ref{m:1}), we obtain
%
\begin{equation}\label{Reh(t)}
\mathfrak{Re} [h(t)]=\Biggl[
\prod_{j=1}^3\frac{d_j(t,k)}{\sqrt{1+4\Delta_j^2t^2}}\Biggr]\cos
\Biggl( \sum_{j=1}^3\bigl( \varphi_j (t)+k^2t\psi_j (t)\bigr)\Biggr),
\end{equation}
where, for $j=1,2,3$:
\begin{eqnarray*}
&\bullet&\hspace*{3pt} d_j(t,k)=\exp\biggl[-\frac{k^2t^2}{2}\frac
{(s_{1j}+s_{2j})^2}{1+4\Delta
_j^2t^2} \biggr];
\\
&\bullet&\hspace*{3pt}
\varphi_j (t)=\frac{1}{2}\arctan(2\Delta_j t), \qquad 0<\varphi_j<\uppi/4;
\\
&\bullet&\hspace*{3pt}
\psi_j (t)=\frac{1}{3}-t^2\frac{(s_{1j}+s_{2j})^2\Delta
_j}{1+4\Delta_j^2t^2}.
\end{eqnarray*}
%
%
Introducing these expressions into (\ref{f:d2}) and using
(\ref{f:d0}), we obtain a new formula which has the form of a rather
complicated integral. However, it is well adapted to numerical
evaluation. On the other hand, this formula allows us to compute the equivalent
as $k\rightarrow0$ of the expectation of the total number of
specular points under the Longuet-Higgins approximation. In fact, a
first-order expansion of the terms in the integrand gives a somewhat
more accurate result, one that we now state as a theorem.

\begin{theo}\label{expspecdim2}
%
\begin{equation}\label{explonguetdim2}
\mathbb{E}(\mathit{SP}_2(\mathbb{R}^2)) = \frac{m_2}{k^2}+\mathrm{O}(1),
\end{equation}
where
%
\begin{eqnarray}\label{2dimlonghigg}
m_2&=&\int_0^{+\infty}\frac{1-[\prod_{j=1}^3(1+4\Delta
_j^2t^2)]^{-1/2}
\cos(\sum_{j=1}^3\varphi_j(t) )}{t^2}\,\mathrm{d}t\nonumber
\\
&=&\int_0^{+\infty}\frac{1-2^{-3/2} [\prod_{j=1}^3
(A_j\sqrt{1+A_j} )]( 1-B_1B_2-B_2B_3-B_3B_1)
}{t^2}\,\mathrm{d}t,
\\
A_j&=&A_j(t)=(1+4\Delta_j^2t^2
)^{-1/2},\qquad B_j=B_j(t)=\sqrt{(1-A_j)/(1+A_j)}.\nonumber
\end{eqnarray}
\end{theo}

Note that $m_2$ depends only on the eigenvalues
$\Delta_1,\Delta_2,\Delta_3$ and is easily computed numerically. We
have performed a numerical computation using a standard sea model with
a Jonswap spectrum and spread function $\cos(2 \theta)$. It
corresponds to the default parameters of the Jonswap function of the
toolbox WAFO \cite{wafo}. The variance matrix of the gradient and the
matrix $\Sigma$ are, respectively,
\[
10^{ -4}
\pmatrix{
114 &0\cr
0& 81\cr
}
,\qquad \Sigma=10^{ -4}
\pmatrix{
9& 3& 0 \cr
3& 11& 0 \cr
0& 0 & 3 \cr
}
.
\]

\begin{figure}[t]

\includegraphics{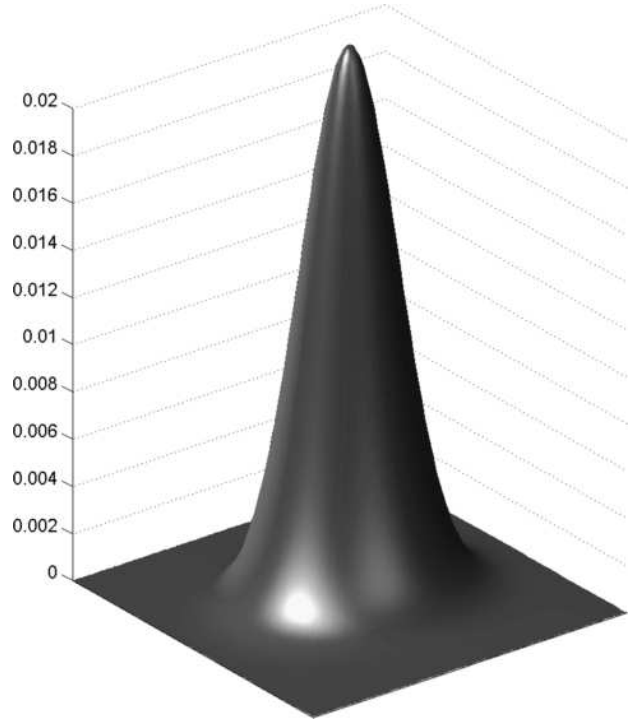}

\caption{Intensity function of the specular points for the
Jonswap spectrum.}\label{f:2}
\end{figure}

The integrand in (\ref{f:spec22}) is displayed in Figure \ref{f:2} as
a function of the two space variables $x,y$. The value of the
asymptotic parameter $m_2$ is $2.527 10 ^{-3}$.

We now consider the variance of the total number of
specular points in two dimensions, looking for analogous results to
the one-dimensional case (i.e., Theorem \ref{asymptoticvar}), in view
of their interest for
statistical applications. It turns out that the computations become
much more complicated. The statements on variance and speed of
convergence to zero of the coefficient of variation that we give
below include only the order of the asymptotic behavior in the
Longuet-Higgins approximation, but not the constant. However, we
still consider them to be useful. If one refines the computations,
rough bounds can be given on the generic constants in Theorem
\ref{varspecdim2} on the basis
of additional
hypotheses on the random field.

We assume that the real-valued, centered, Gaussian stationary random
field $\{W(\mathbf{x})\dvtx \mathbf{x}\in\mathbb{R}^2\}$ has paths of class
$C^3$, the distribution of $W'(\mathbf{0})$ does not degenerate
(i.e., $\operatorname{Var}(W'(\mathbf{0}))$ is invertible). Moreover,
let us
consider $W''(\mathbf{0})$, expressed in the reference system $xOy$
of $\mathbb{R}^2 $ as the $2\times2$ symmetric centered Gaussian random
matrix
\[
W''(\mathbf{0})=
\pmatrix{
W_{xx}(\mathbf{0}) & W_{xy}(\mathbf{0}) \cr
W_{xy}(\mathbf{0}) & W_{yy}(\mathbf{0}) \cr
}
.
\]
The function
\[
\mathbf{z}\rightsquigarrow\Delta(\mathbf{z})=\det[
\operatorname{Var}(
W''(0)\mathbf{z} ) ],
\]
defined on $\mathbf{z}=(z_1,z_2)^\mathrm{T} \in\mathbb{R}^2$, is a non-negative
homogeneous polynomial of degree $4$ in the pair $z_1,z_2$. We will
assume the non-degeneracy condition
%
\begin{equation}\label{nondegdelta}
\min\{\Delta(\mathbf{z})\dvtx \|\mathbf{z}\|=1 \}=\underline{\Delta
}>0.
\end{equation}

\begin{theo}\label{varspecdim2}
Let us assume that $\{W(\mathbf{x})\dvtx \mathbf{x}\in\mathbb{R}^2 \}$ satisfies
the above conditions and that it is also $\delta$-dependent, $\delta
>0 $, that is, $\mathbb{E}(W(\mathbf{x})W(\mathbf{y}))=0$
whenever $\|\mathbf{x}-\mathbf{y}\|>\delta.$ Then, for $k$ small enough,
%
\begin{eqnarray}\label{var}
\operatorname{Var}(\mathit{SP}_2(\mathbb{R}^2) )\leq \frac{L}{k^2},
\end{eqnarray}
where $L$ is a positive constant depending on the law of the
random field.

Moreover, for $k$ small enough, by using the result of Theorem \ref
{expspecdim2} and (\ref{var}), we get
\[
\frac{\sqrt{\operatorname{Var}(\mathit{SP}_2(\mathbb{R}^2)
)}}{\mathbb{E}(\mathit{SP}_2(\mathbb{R}^2)
)}\leq L_1k,
\]
where $L_1$ is a new positive constant.
\end{theo}

\begin{pf} To simplify notation, let
us denote $T=\mathit{SP}_2( \mathbb{R}^2)$. We have
%
\begin{equation}\label{vardim2}
\operatorname{Var}(T)=\mathbb{E}\bigl(T(T-1)\bigr)+\mathbb{E}(T)-[\mathbb{E}(T)]^2.
\end{equation}

We have already computed the equivalents as $k\rightarrow0$ of the
second and third term in the right-hand side of (\ref{vardim2}). Our
task in what follows is to consider the first term.

The proof is performed along the same lines as the one of Theorem
\ref{asymptoticvar}, but instead of applying a Rice formula for the
second factorial moment of the number of crossings of a
one-parameter random process, we need \cite{aw}, Theorem 6.3, for the
factorial moments of a 2-parameter random field.
We have
\begin{eqnarray*}\label{2dim2}
\mathbb{E}\bigl(T(T-1)\bigr)
&=&\int\hspace*{-2pt}\int_{\mathbb{R}^2 \times
\mathbb{R}^2}\mathbb{E}\bigl(|\det\mathbf{Y}'(\mathbf{x})||\det
\mathbf
{Y}'(\mathbf{y})||
\mathbf{Y}(\mathbf{x})=\mathbf{0},\mathbf{Y}(\mathbf{y})=\mathbf{0}
\bigr)
\\
&&{}\hspace*{36pt}\times p_{\mathbf{Y}(\mathbf{x}),\mathbf{Y}(\mathbf{y})}(\mathbf
{0},\mathbf{0})\,\mathrm{d}\mathbf{x}\,\mathrm{d}\mathbf{y}\nonumber
\\
&=&\int\hspace*{-2pt}\int_{\|\mathbf{x}-\mathbf{y}\|>\delta}\cdots \,\mathrm{d}\mathbf
{x}\,\mathrm{d}\mathbf{y}+\int\hspace*{-2pt}\int_{\|\mathbf{x}-\mathbf{y}\|\leq
\delta}\cdots
\,\mathrm{d}\mathbf{x}\,\mathrm{d}\mathbf{y}=J_1+J_2.
\end{eqnarray*}

For $J_1$, we proceed as in the proof of Theorem \ref{asymptoticvar},
using the $\delta$-dependence and the evaluations leading to the
statement of Theorem \ref{expspecdim2}. We obtain
%
\begin{equation}\label{j1}
J_1=\frac{m_2^2}{k^4}+\frac{\mathrm{O}(1)}{k^2}.
\end{equation}
One can show that under the hypotheses of the theorem, for small $k$,
one has
%
\begin{equation}\label{j2}
J_2=\frac{\mathrm{O}(1)}{k^2}.
\end{equation}
We refer the reader to \cite{azlomw} for the lengthy computations
leading to this inequality. In view of (\ref{vardim2}), (\ref
{explonguetdim2}) and
(\ref{j1}), this suffices to prove the theorem.
\end{pf}

\section{The distribution of the normal to the level curve} \label{s:normal}

Let us consider a modeling of the sea $ W(x,y,t)$ as a function
of two space variables and one time variable. Usual models are centered
Gaussian stationary with a particular form of the spectral measure $\mu$
that is presented, for example, in \cite{aw}. We denote the covariance
by $ \Gamma(x,y,t)=\mathbb{E}(W(0,0,0)W(x,y,t))$.

In practice, one is frequently confronted with the following situation:
several pictures of the sea on time over an interval $[0,T]$ are
stocked and some properties or magnitudes are observed. If the time
$T$ and the number of pictures are large, and if the process is
ergodic in time, then the frequency of pictures that satisfy a certain
property will converge to the probability of this property happening
at a fixed time.

Let us illustrate this with the angle of the normal to the level
curve at a point ``chosen at random''. We first consider the number
of crossings of a level $u$ by the process $W(\cdot,y,t)$ for fixed
$t$ and $y$, defined as
\[
N_{[0,M_1]}^{W(\cdot, y,t)}(u)=\#\{ x \dvtx 0\le x \le
M_1; W(x,y,t)=u\}.
\]
We are interested in computing the total number of crossings per
unit time when integrating over $ y \in[0,M_2]$, that is,
%
\begin{equation}\label{hori}
\frac1T \int_0^T \mathrm{d}t \int_0^{M_2}N_{[0,M_1]}^{W(\cdot,y,t)}(u)\,\mathrm{d}y.
\end{equation}
If the ergodicity assumption in time holds true, then we can conclude
that a.s.
\[
\frac1T
\int_0^T \mathrm{d}t \int_0^{M_2}N_{[0,M_1]}^{W(\cdot,y,t)}(u)\,\mathrm{d}y \to M_1
\mathbb{E}\bigl(N_{[0,M_1]}^{W(\cdot,0,0)}(u)\bigr)
=\frac{M_1M_2}{\uppi
}\sqrt{\frac{\lambda_{200}}{\lambda_{000}}}\mathrm{e}^{-1/2
{u^2/\lambda_{000}}},
\]
where
\[
\lambda_{\mathit{abc}} = \int_{\mathbb{R}^3} \lambda_x^a\lambda_y^b\lambda_t^c
\,\mathrm{d}\mu(\lambda_x,\lambda_y,\lambda_t)
\]
are the spectral moments of $W$. Hence, on the basis of the quantity
(\ref{hori}), for large $T$, one
can make inference about the value of certain parameters of the law
of the random field. In this example, these are the spectral moments
$\lambda_{200}$ and $\lambda_{000}$.

If two-dimensional level information is available, one can work
differently because there exists an interesting relationship with Rice
formulae for level curves that we
explain in what follows. We can write ($\mathbf{x}=(x,y)$)
\[
W'(\mathbf{x},t)=\|W'(\mathbf{x},t)\|(\cos\Theta(\mathbf
{x},t),\sin
\Theta(\mathbf{x},t))^T.
\]
Using a Rice formula, more precisely, under conditions of \cite{aw},
Theorem 6.10,
%
\begin{equation}\label{f:chichi}
\mathbb{E}\biggl[
\int_0^{M_2}N_{[0,M_1]}^{W(\cdot,y,0)}(u)\,\mathrm{d}y\biggr]=\mathbb{E}
\biggl[\int
_{\mathcal{
C}_Q(0,u)}|\cos\Theta(\mathbf{x},0)|\,\mathrm{d}\sigma_1\biggr]
=\frac{\sigma_2(Q)}{\uppi}\sqrt{\frac{\lambda_{200}}{\lambda
_{000}}} \mathrm{e}^{-{u^2/(2\lambda_{000})}},\
\end{equation}
where $ Q = [0,M_1]\times[0,M_2] $. We have a similar formula when
we consider sections of the set $ [0,M_1]\times[0,M_2] $ in the
other direction. In fact, (\ref{f:chichi}) can be generalized to
obtain the Palm distribution of the angle $ \Theta$.

\begin{figure}[t]

\includegraphics{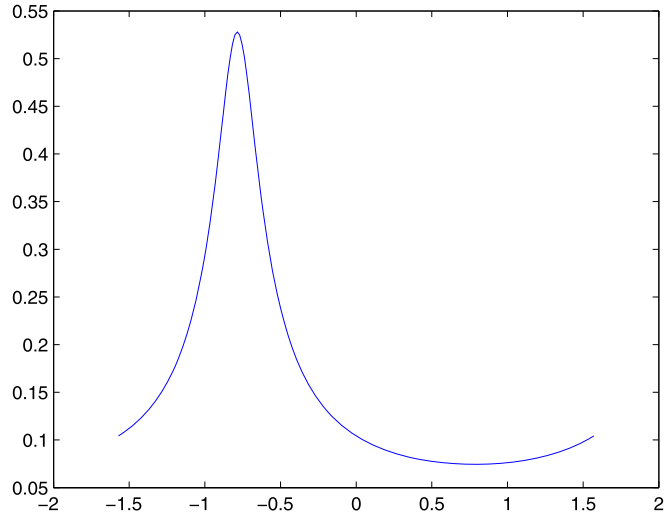}

\caption{Density of the Palm distribution of the angle of the normal
to the level curve in the case $\gamma=0.5$ and
$\kappa=\uppi/4$.} \label{fi:3}
\end{figure}

Set $ h_{\theta_1,\theta_2}=\mathbb{I}_{[
\theta_1 , \theta_2]}$ and, for
$-\uppi\le\theta_1<\theta_2\le\uppi$, define
%
\begin{eqnarray}\label{f:hh}
F(\theta_2)-F(\theta_1)&:=& \mathbb{E}\bigl(\sigma_1\bigl(\{\mathbf
{x}\in Q\dvtx
W(\mathbf{x},0)=u  ;  \theta_1\le\Theta(\mathbf{x},s)\le
\theta_2\}\bigr)\bigr) \nonumber
\\
&\phantom{:}\!=& \mathbb{E}\biggl(\int_{\mathcal{C}_Q(u,s)}h_{\theta_1,\theta_2}(
\Theta(\mathbf{x},s))\,\mathrm{d}\sigma_1(\mathbf{x})\,\mathrm{d}s\biggr)
\\
&\phantom{:}\!=&\sigma_2(Q)\mathbb{E}\biggl[h_{\theta_1,\theta_2}\biggl(\frac{\partial_{y}W}
{\partial_{x}W}\biggr)\bigl((\partial_{x}W)^2+(\partial_{y}W)^2\bigr)^{1/2}\biggr] \frac
{\exp(-{u^2/(2\lambda_{00})})}{\sqrt{
2\uppi\lambda_{000}}}.\nonumber
\end{eqnarray}
Defining $\Delta=\lambda_{200}\lambda_{020}-\lambda_{110}$ and
assuming $\sigma_2(Q)=1$ for ease of notation, we readily obtain
\begin{eqnarray*}
&&F(\theta_2)-F(\theta_1)
\\
&&\quad=\frac{ \mathrm{e}^{-{u^2/(2\lambda_{000})}}}
{(2\uppi)^{3/2}(\Delta)^{1/2}\sqrt{\lambda_{000}}}\int_{\mathbb{R}
^2}h_{\theta_1,\theta_2}(\Theta)
\sqrt{x^2+y^2}
\mathrm{e}^{-(1/({2\Delta}))(\lambda_{02}x^2-2\lambda
_{11}xy+\lambda_{20}y^2)}\,\mathrm{d}x\,\mathrm{d}y
\\
&&\quad=\frac{ \mathrm{e}^{-{u^2/(2\lambda_{00})}}}{(2\uppi)^{3/2}(\lambda
_{+}\lambda_{-})^{1/2}\sqrt{\lambda_{000}}}
\\
&&{}\qquad \times\int_0^{+\infty}\hspace*{-2pt}\int_{\theta_1}^{\theta_2}\rho^2
\exp\biggl(-\frac{\rho^2}{2\lambda_{+}\lambda_{-}}\bigl(\lambda_{+}\cos
^2(\varphi-\kappa)+\lambda_{-}
\sin^2(\varphi-\kappa)\bigr)\biggr)\,\mathrm{d}\rho \,\mathrm{d}\varphi,
\end{eqnarray*}
where $\lambda_{-}\le\lambda_{+}$ are the eigenvalues of the
covariance matrix of the random vector
$(\partial_{x}W(0,0,0),\partial_{y}W(0,0,0))$ and $\kappa$ is
the angle of the eigenvector associated with $ \gamma^+$. Noting
that the exponent in the integrand can be written as $
1/\lambda_-(1 -\gamma^2 \sin^2( \varphi-\kappa)) $ with $\gamma^2 :=
1- \lambda_+/\lambda_-$ and that
\[
\int_0^{+\infty} \rho^2 \exp\biggl( -\frac{H \rho^2}{2}\biggr) =
\sqrt{\frac{\uppi}{2H}},
\]
it is easy to obtain that
\[
F(\theta_2)-F(\theta_1) = (\mathit{const}) \int_{\theta_1}^{\theta_2} \bigl(1
-\gamma^2 \sin^2( \varphi-\kappa)\bigr)^{-1/2} \,\mathrm{d} \varphi.
\]
From this relation, we get the density $g(\varphi)$ of the Palm
distribution, simply by dividing by the total mass:
%
\begin{equation}\label{f:palm}
g(\varphi) = \frac{(1 -\gamma^2 \sin^2( \varphi
-\kappa))^{-1/2}}{ \int_{-\uppi}^{\uppi} (1 -\gamma^2 \sin^2(
\varphi-\kappa))^{-1/2}\,\mathrm{d} \varphi} =\frac{(1 -\gamma^2
\sin^2( \varphi-\kappa))^{-1/2}}{ 4 \mathcal{K}(\gamma^2)} .
\end{equation}
Here, $ \mathcal{K}$ is the complete elliptic integral of the first
kind. This density characterizes the distribution of the angle of
the normal at a point chosen ``at random'' on the level curve.
In the case of a random field which is isotropic in $(x,y)$, we have
$\lambda_{200}=\lambda_{020}$ and, moreover, $\lambda_{110}=0$, so
that $g$ turns out to be the uniform density over the circle
(Longuet-Higgins says that over the contour, the ``distribution''
of the angle is uniform (cf.~\cite{lo2lo2}, page 348)).
We have performed the numerical computation of the density (\ref
{f:palm}) for an
anisotropic process with $\gamma=0.5$, $\kappa=\uppi/4$.
Figure \ref{fi:3} displays the densities of the Palm distribution
of the angle showing a large departure from the uniform
distribution.

Let us turn to ergodicity. For a given subset $Q$ of $\mathbb{R}^2$
and each
$t$, let us define
$\mathcal{A}_t=\sigma\{W(x,y,t)\dvtx\tau>t ;(x,y)\in
Q\}$ and consider the $\sigma$-algebra of $t$-invariant events $
\mathcal{A}=\bigcap\mathcal{A}_t $. We assume that for each pair
$(x,y)$, $\Gamma(x,y,t)\rightarrow0$ as $t\rightarrow+\infty.$ It
is well known that under this condition, the $\sigma$-algebra $
\mathcal{A} $ is trivial, that is, it only contains events having
probability zero or one (see, e.g., \cite{crle}, Chapter 7). This has
the following important consequence in our context. Assume that the set
$Q$ has a smooth boundary and, for simplicity, unit Lebesgue measure.
Let us consider
%
\begin{equation}\label{e:z1}
Z(t)=\int_{\mathcal{C}_Q(u,t)}H(\mathbf{x},t)\,\mathrm{d}\sigma
_1(\mathbf{x})
\end{equation}
with $ H(\mathbf{x},t)=\mathcal{H}( W(\mathbf
{x},t),\nabla W (\mathbf{x},t ))$, where $\nabla W =(W_x,W_y)$
denotes the gradient in the space variables and $\mathcal{H} $ is some
measurable function such that the integral is well defined. This is
exactly our case in
(\ref{f:hh}). The process $\{Z(t)\dvtx t\in\mathbb{R}\}$
is strictly stationary and, in our case, has a finite mean and is
Riemann-integrable. By the Birkhoff--Khintchine ergodic theorem
(\cite{crle}, page 151), a.s.~as $T\rightarrow+\infty$,
\[
\frac1T\int_0^T Z(s)\,\mathrm{d}s\to\mathbb{E}_{\mathcal{B}}[Z(0)],
\]
where $ \mathcal{B}$ is the $\sigma$-algebra of $t$-invariant
events associated with the process $Z(t)$. Since for each~$t$, $Z(t)$
is $\mathcal{A}_t$-measurable, it follows that $ \mathcal{B} \subset
\mathcal{A}$ so that $\mathbb{E}_{\mathcal{B}}[Z(0)]=\mathbb
{E}[Z(0)] $.
On the other hand, the Rice formula yields
(taking into account the fact that stationarity of $\mathcal{W} $
implies that $W(\mathbf{0},0)$ and $\nabla W(\mathbf{0},0)$ are independent)
\begin{eqnarray*}
\mathbb{E}[Z(0)]=
\mathbb{E}[\mathcal{H}( u,\nabla W(\mathbf{0},0))\|\nabla
W(\mathbf{0},0)\|]p_{W(\mathbf{0},0)}(u).
\end{eqnarray*}
We consider now the central limit theorem. Let us define
%
\begin{equation}\label{e:z2}
\mathcal{Z}(t)=\frac1t\int_0^t [ Z(s) -\mathbb{E}( Z(0))
]\, \mathrm{d}s.
\end{equation}
To compute the variance of $\mathcal{Z}(t)$, one can again use the
Rice formula for the first moment of integrals over level sets, this
time applied to the $\mathbb{R}^2$-valued random field with parameter
in $\mathbb{R}
^4$, $\{(W(\mathbf{x}_1,s_1),W(\mathbf{x}_2,s_2))^\mathrm{T}\dvtx (\mathbf
{x}_1,\mathbf{x}_2)\in Q\times Q,s_1,s_2 \in[0,t] \}$ at the level $(u,u)$.
We get
\[
\operatorname{Var}\mathcal{Z}(t)=\frac2t\int_0^t
\biggl(1-\frac st\biggr)I(u,s)\,\mathrm{d}s,
\]
where
\begin{eqnarray*}
I(u,s)&=&\int_{Q^2} \mathbb{E}[H(\mathbf{x}_1, 0) H(\mathbf
{x}_2, s) \|
\nabla W(\mathbf{x}_1,0)\|
\| \nabla W(\mathbf{x}_2,s)\|
|W(\mathbf{x}_1,0)=u ;W(\mathbf{x}_2,s)=u ]
\\
&&{}\hspace*{14pt}\times p_{W(\mathbf{x}_1,0),W(\mathbf{x}_2,s)}(u,u)\, \mathrm{d}\mathbf{x}_1
\,\mathrm{d}\mathbf{x}_2
- \bigl(\mathbb{E}[\mathcal{H}( u,\nabla
W(\mathbf
{0},0))\|\nabla W(\mathbf{0},0)\|]p_{W(\mathbf{0},0)}(u) \bigr)^2.
\end{eqnarray*}
Assuming that the given random field is time-$\delta$-dependent, that
is,
$\Gamma(x,y,t)=0\, \forall (x,y) $ whenever $t>\delta$, we readily
obtain
%
\begin{equation}\label{e:sig}
t\operatorname{Var}\mathcal{\mathcal{Z}}(t)\to2\int_0^{\delta}
I(u,s)\,\mathrm{d}s:=\sigma^2(u) \qquad\mbox{as }t\to\infty.
\end{equation}
Now, using a variant of the Hoeffding--Robbins theorem \cite{hoero}
for sums of
$\delta$-dependent random variables, we can establish the following theorem.

\begin{theo}
Assume that the random field $W$ and the function $H$
satisfy the conditions of \cite{aw}, Theorem
6.10. Assume, for simplicity, that $Q$ has Lebesgue measure.
Then:
\begin{longlist}[(ii)]
\item[(i)] if the covariance $\gamma(x,y,t)$ tends to zero as $t \to
+\infty$ for every value of $(x,y) \in Q$, we have
\[
\frac1T\int_0^T Z(s)\,\mathrm{d}s\to\mathbb{E}[\mathcal{H}( u,\nabla
W(\mathbf
{0},0))\|\nabla W(\mathbf{0},0)\|]p_{W(\mathbf{0},0)}(u),
\]
where $ Z(t)$ is defined by (\ref{e:z1}).

\item[(ii)] if the random field $W$ is $\delta$-dependent in the sense
above, we have
\[
\sqrt t \mathcal{Z}(t)\quad\Longrightarrow\quad N(0,\sigma^2(u)),
\]
where $\mathcal{Z}(t)$ is defined by (\ref{e:z2}) and $\sigma^2(u)$
by (\ref{e:sig}).
\end{longlist}
\end{theo}

\section{Application to dislocations of wavefronts}\label{section4}

In this section, we follow the article
\cite{be1be1} by Berry and Dennis. Dislocations are lines in space or
points in the plane where
the phase $\chi$ of the complex scalar wave
$\displaystyle\psi(\mathbf{x},t)=\rho(\mathbf{x},t)e^{\mathrm{i}\chi
(\mathbf{x},t)}$
is undefined. With respect to light, they are lines of darkness;
with respect to sound, threads of silence. Here, we only consider
two-dimensional space variables $\mathbf{x} = (x_1,x_2)$.

It is convenient to express $\psi$ by means of its real and
imaginary parts:
\[
\psi(\mathbf{x},t)=\xi(\mathbf{x},t)+\mathrm{i}\eta(\mathbf{x},t).
\]
Thus,
dislocations are the intersection of the surfaces $
\xi(\mathbf{x},t)=0$ and $\eta(\mathbf{x},t)=0.$

Let us quote the authors of \cite{be1be1}: ``Interest in optical
dislocations has recently revived, largely as a result of experiments
with laser fields. In low-temperature
physics, $\psi(\mathbf{x},t)$ could represent the complex order
parameter associated with
quantum flux lines in a superconductor or quantized vortices in a
superfluid'' (cf.~\cite{be1be1} and the references therein).

In what follows, we assume an isotropic Gaussian model. This means that
we will
consider the wavefront as an isotropic Gaussian field
\[
\psi(\mathbf{x},t)=\int_{{\mathbb{R}}^2}\exp{(\mathrm{i}[ \langle\mathbf{k}\cdot
\mathbf{x}
\rangle-c
|\mathbf{k}|t])}\biggl(\frac{\Pi(|\mathbf{k}|)}{|\mathbf{k}|}\biggr)^{1/2}\,\mathrm{d}W(\mathbf{k}),
\]
where
$\mathbf{k}=(k_1,k_2)$, $|\mathbf{k}|=\sqrt{k_1^2+k_2^2}$, $\Pi(k)$ is the
isotropic spectral density and\vspace*{-2pt} $W=(W_1+\mathrm{i}W_2)$ is a standard complex
orthogonal Gaussian measure on $\mathbb{R}^2$ with unit variance.
We are only interested in $t=0$ and we put
$\xi(\mathbf{x}):=\xi({\mathbf
x},0)$ and $\eta(\mathbf{x}):=\eta(\mathbf{x},0)$. We have, setting $k =
|\mathbf{k}|$,
%
\begin{eqnarray}\label{process1}
\xi(\mathbf{x})&=&\int_{{\mathbb{R}}^2}\cos( \langle\mathbf{k}\cdot{\mathbf
x}\rangle)\biggl(\frac{\Pi(k)}{k}\biggr)^{1/2}\,\mathrm{d}W_1(\mathbf{k})-\int_{{\mathbb
{R}}^2} \sin
(\langle{\mathbf
k}\cdot\mathbf{x}\rangle)\biggl(\frac{\Pi(k)}{k}\biggr)^{1/2}\,\mathrm{d}W_2(\mathbf{k}),
\\\label{process2}
\eta(\mathbf{x})&=&\int_{{\mathbb{R}}^2}\cos(\langle\mathbf{k}\cdot{\mathbf
x}\rangle)\biggl(\frac{\Pi(k)}{k}\biggr)^{1/2}\,\mathrm{d}W_2(\mathbf{k})+\int_{{\mathbb
{R}}^2} \sin
(\langle{\mathbf
k}\cdot\mathbf{x}\rangle)\biggl(\frac{\Pi(k)}{k}\biggr)^{1/2}\,\mathrm{d}W_1(\mathbf{k}).
\end{eqnarray}
The covariances are
%
\begin{equation}\label{covariance}
\mathbb{E} [\xi(\mathbf{x})\xi(\mathbf{x}')]=\mathbb{E} [\eta({\mathbf
x})\eta(\mathbf{x}')]
= \rho(|\mathbf{x}-\mathbf{x}'|) :=
\int_0^\infty J_0(k|\mathbf{x}-\mathbf{x}'|)\Pi(k)\,\mathrm{d}k,
\end{equation}
where $J_{\nu}(x)$ is the Bessel function of the first kind of order
$\nu$. Moreover, $\mathbb{E} [\xi(\mathbf{r}_1)\eta(\mathbf{r}_2)]=0.$

\subsection{Mean number of dislocation points}

Let us denote by
$\{\mathbf{Z}(\mathbf{x})\dvtx \mathbf{x} \in\mathbb{R}^2\}$ a random field
having values in $\mathbb{R}^2$, with coordinates $\xi(\mathbf{x}),
\eta(\mathbf{x})$, which are two independent Gaussian stationary
isotropic random fields with the same distribution. We are
interested in the expectation of the number of dislocation points
\[
d_2 := \mathbb{E}[ \#\{\mathbf{x} \in S \dvtx\xi(\mathbf{x})
=\eta(\mathbf{x})= 0\}],
\]
where $S$ is a subset of the parameter space having area equal to
$1$.

Without loss of generality, we may assume that $
\operatorname{Var}(\xi(\mathbf{x})) = \operatorname{Var}(\eta
(\mathbf{x})) =1 $ and for the
derivatives, we set $ \lambda_2= \operatorname{Var}(\eta_i(\mathbf
{x}))=\operatorname{Var}
(\xi_i(\mathbf{x}))$, $ i=1,2$. Then, according to the
Rice formula,
\[
d_2 =\mathbb{E}[ |\det(\mathbf{Z}'(\mathbf{x})) |/ \mathbf
{Z}(\mathbf{x})
=0] p_{\mathbf{Z}(\mathbf{x})}(0).
\]
An easy Gaussian computation gives $d_2 = \lambda_2 /(2\uppi) $
(\cite{be1be1}, formula (4.6)).

\subsection{Variance}

Again, let
$S$ be a measurable subset of $\mathbb{R}^2$ having Lebesgue measure
equal to $1$. We have
\[
\operatorname{Var}( N^\mathbf{Z}_S(\mathbf{0}))=\mathbb{E}
\bigl(N^\mathbf{Z}_S(\mathbf{0})\bigl(N^\mathbf{Z}_S(\mathbf{0})-1
\bigr) \bigr)+d_2-d_2^2
\]
and for the first term, we use the Rice formula for the second
factorial moment (\cite{aw}, Theorem~6.3), that is,
\[
\mathbb{E}\bigl(N^\mathbf{Z}_S(\mathbf{0})\bigl(N^\mathbf
{Z}_S(\mathbf{0})-1
\bigr) \bigr) = \int_{S^2}
A(\mathbf{s}_1,\mathbf{s}_2)\, \mathrm{d}\mathbf{s}_1\, \mathrm{d}\mathbf{s}_2,
\]
where
\[
A(\mathbf{s}_1,\mathbf{s}_2) =\mathbb{E}[ |\det\mathbf
{Z}'(\mathbf
{s}_1) \det\mathbf{Z}'(\mathbf{s}_2)| |
\mathbf{Z}(\mathbf{s}_1) = \mathbf{Z}(\mathbf{s}_2)=\mathbf{0}_2
] p_{\mathbf{
Z}(\mathbf{s}_1), \mathbf{Z}(\mathbf{s}_2)}(\mathbf{0}_4).
\]
Here, $ \mathbf{0}_p$ denotes the null vector in dimension $p$.

Taking into account the fact that the law of the random field $\mathbf
{Z}$ is invariant
under translations and orthogonal transformations of $\mathbb{R}^2$,
we have
\[
A(\mathbf{s}_1,\mathbf{s}_2) = A((0,0),(r,0) ) = A(r) \qquad
\mbox{with } r= \| \mathbf{s}_1-\mathbf{s}_2\|.
\]

The Rice function $ A(r)$ has two intuitive interpretations. First, it
can be viewed as
\[
A(r)= \lim_{\epsilon\to0} \frac{1}{ \uppi^2
\epsilon^4}\mathbb{E}[N( B( (0,0), \epsilon))\times
N( B(
(r,0), \epsilon))].
\]
Second, it is the density of the
Palm distribution, a generalization of the horizontal window conditioning
of the number of zeros of $\mathbf{Z}$ per unit
surface, locally around the point $(r,0)$, conditionally on the
existence of a
zero at $(0,0)$ (see \cite{crle}). In \cite{be1be1}, $A(r)/d_2^2$
is called the ``correlation function''.

To compute $A(r)$, we denote by $\xi_1,\xi_2,\eta_1, \eta_2$
the partial derivatives of $\xi, \eta$ with respect to first and
second coordinate. Therefore,
%
\begin{eqnarray}\label{corre}
A(r)&=& \mathbb{E}[ |\det\mathbf{Z}'(0,0) \det\mathbf
{Z}'(r,0)| |
\mathbf{Z}(0,0) = \mathbf{Z}(r,0)=\mathbf{0}_2 ] p_{\mathbf{
Z}(0,0), \mathbf{Z}(r,0)}(\mathbf{0}_4)\nonumber
\\
&=&\mathbb{E}[ |(\xi_1\eta_2-\xi_2\eta_1
)(0,0)(
\xi_1\eta_2-\xi_2\eta_1)(r,0) | |\mathbf
{Z}(0,0) =
\mathbf{Z}(r,0)=\mathbf{0}_2]
\\
&&{}\times  p_{\mathbf{Z}(0,0),
\mathbf{Z}(r,0)}(\mathbf{0}_4).\nonumber
\end{eqnarray}
The density is easy to compute:
\[
p_{ \mathbf{Z}(0,0), \mathbf{Z}(r,0)}(\mathbf{0}_4) = \frac
{1}{(2\uppi)^2(1-\rho^2(r))},
\qquad\mbox{where } \rho(r)=\int_0^\infty J_0(kr)\Pi(k)\,\mathrm{d}k.
\]
The conditional expectation turns out to be more difficult to
calculate, requiring a long computation (we again refer to \cite
{azlomw} for the details). We obtain the following formula (that can
be easily compared to the formula in \cite{be1be1} since we are using
the same notation):
\[
A(r)= \frac{A_1}{4 \uppi^3(1-C^2)} \int_{-\infty}^{\infty} \frac
1{t^2}\biggl[1-\frac1{(1+t^2)}
\frac{(Z_2-2Z_1^2t^2)}{ Z_2\sqrt{(Z_2-Z_1^2t^2)}}\biggr]\,\mathrm{d}t,
\]
\eject
where we have defined
\begin{eqnarray*}
C &:=& \rho(r), \qquad E = \rho'(r), \qquad H= -E/r, \qquad F =-\rho''(r),\qquad
F_0=-\rho''(0),
\\
 A_1&=&F_0\biggl(F_0 -\frac{E^2}{1-C^2}\biggr),\qquad A_2 = \frac{ H}{F_0} \frac{ F(
1-C^2) - E^2C}{F_0( 1-C^2) -E^2 },
\\
Z &=& \frac{F_0^2-H^2}{F_0^2} \biggl[ 1-\biggl(F -\frac{E^2
C}{1-C^2}\biggr)^2\cdot
\biggl(F_0 -\frac{E^2}{1-C^2}\biggr)^{-2}
\biggr],
\\
Z_1&=&\frac{A_2}{1+Zt^2},\qquad Z_2=\frac{1+t^2}{1+Zt^2}.
\end{eqnarray*}

\section*{Acknowledgement} This work has received financial support from
the European Marie Curie Network SEAMOCS.

\printhistory

\end{document}